\documentclass[11pt]{article}
\usepackage[letterpaper, margin=1.1in]{geometry}
\usepackage[utf8]{inputenc}
\usepackage{amsmath,amsthm,amssymb,amsopn}
\usepackage{authblk}
\usepackage{lipsum}

\newcommand\blfootnote[1]{%
  \begingroup
  \renewcommand\thefootnote{}\footnote{#1}%
  \addtocounter{footnote}{-1}%
  \endgroup
}

\newbox\gnBoxA
\newdimen\gnCornerHgt
\setbox\gnBoxA=\hbox{$\ulcorner$}
\global\gnCornerHgt=\ht\gnBoxA
\newdimen\gnArgHgt
\def\gnmb #1{\setbox\gnBoxA=\hbox{$#1$}\gnArgHgt=\ht\gnBoxA\ifnum     \gnArgHgt<\gnCornerHgt \gnArgHgt=0pt\else \advance \gnArgHgt by -\gnCornerHgt\fi \raise\gnArgHgt\hbox{$\ulcorner$} \box\gnBoxA \raise\gnArgHgt\hbox{$\urcorner$}}
\theoremstyle{plain}

\theoremstyle{remark}

\theoremstyle{definition}

\begin{document}

\title{Truth, Disjunction, and Induction}
\author[*]{Ali~Enayat}
\author[**]{Fedor~Pakhomov}
\affil[*]{\small {Department of Philosophy, Linguistics, and the Theory of Science, University of Gothenburg, Gothenburg, Sweden~~\texttt{ali.enayat@gu.se}}}

\affil[**]{\small{Steklov Mathematical Institute of Russian Academy of Sciences, Moscow, Russia~~\texttt{pakhfn@mi.ras.ru}}}

\renewcommand\Authands{ and }

\maketitle

\begin{abstract}
By a well-known result of Kotlarski, Krajewski, and Lachlan (1981),
first-order Peano arithmetic $\mathsf{PA}$ can be conservatively extended to
the theory $\mathsf{CT}^{-}\mathsf{[PA]}$ of a truth predicate satisfying
compositional axioms, i.e., axioms stating that the truth predicate is
correct on atomic formulae and commutes with all the propositional
connectives and quantifiers. This results motivates the general question of
determining natural axioms concerning the truth predicate that can be added
to $\mathsf{CT}^{-}\mathsf{[PA]}$ while maintaining conservativity over $%
\mathsf{PA}$. Our main result shows that conservativity fails even for the
extension of $\mathsf{CT}^{-}\mathsf{[PA]}$ obtained by the seemingly weak
axiom of disjunctive correctness $\mathsf{DC}$ that asserts that the truth
predicate commutes with disjunctions of arbitrary finite size. In
particular, $\mathsf{CT}^{-}\mathsf{[PA]+DC}$ implies $\mathsf{Con}(\mathsf{%
PA})$.\medskip

Our main result states that the theory $\mathsf{CT}^{-}\mathsf{[PA]+DC}$
coincides with the theory $\mathsf{CT}_{0}\mathsf{[PA]}$ obtained by adding $%
\Delta _{0}$-induction in the language with the truth predicate. This result
strengthens earlier work by Kotlarski (1986) and Cie\'{s}li\'{n}ski (2010).
For our proof we develop a new general form of Visser's theorem on
non-existence of infinite descending chains of truth definitions and prove it by reduction to (L%
\"{o}b's version of) G\"{o}del's second incompleteness theorem, rather than by using the Visser-Yablo paradox, as in Visser's original proof (1989).
\end{abstract}

\blfootnote {Fedor Pakhamov's research is supported
in part by the Young Russian Mathematics award.}



\begin{center}
\textbf{1.~INTRODUCTION}\bigskip
\end{center}

By a theorem of Krajewski, Kotlarski, and Lachlan \cite{KKL}, every
countable recursively saturated model $\mathcal{M}$ of $\mathsf{PA}$ (Peano
Arithmetic) carries a `full satisfaction class', i.e., there is a subset $S$
of the universe $M$ of $\mathcal{M}$ that `decides' the truth/falsity of
each sentence of arithmetic in the sense of $\mathcal{M}$ -- even sentences
of nonstandard length -- while obeying the usual recursive clauses guiding
the behavior of a Tarskian satisfaction predicate. This remarkable theorem
implies that theory $\mathsf{CT}^{-}\mathsf{[PA]}$ (compositional truth over
$\mathsf{PA}$ with induction only for the language $\mathcal{L}_{\mathsf{A}}$
of arithmetic) is \textit{conservative} over $\mathsf{PA}$, i.e., if an $%
\mathcal{L}_{\mathsf{A}}$-sentence $\varphi $ is provable in $\mathsf{CT}^{-}%
\mathsf{[PA]}$, then $\varphi $ is already provable in $\mathsf{PA}$. New
proofs of this conservativity result were given by Visser and Enayat \cite%
{Ali-Albert} using basic model theoretic ideas, and by Leigh \cite{Leigh}
using proof theoretic tools; these new proofs make it clear that in the
Krajewski-Kotlarski-Lachlan theorem the theory $\mathsf{PA}$ can be replaced
by any `base' theory that supports a modicum of coding machinery for
handling elementary syntax.\medskip

On the other hand, it is well-known \cite[Thm.~8.39 and Cor.~8.40]{Halbach
book} that the consistency of $\mathsf{PA}$ (and much more) is readily
provable in the stronger theory $\mathsf{CT[PA]}$, which is the result of
strengthening $\mathsf{CT}^{-}\mathsf{[PA]}$ with the scheme of induction
over natural numbers for all $\mathcal{L}_{\mathsf{A+T}}$-formulae, where $%
\mathcal{L}_{\mathsf{A+T}}:=\mathcal{L}_{\mathsf{A}}\cup \{\mathsf{T}(x)\}.$%
\footnote{$\mathsf{CT[PA]}$ dwarfs $\mathsf{PA}$ in arithmetical strength:
By a classical theorem (discovered by a number of researchers, including
Feferman and Takeuti, and explained in \cite{Halbach book}) the arithmetical
consequences of $\mathsf{CT[PA]}$ are the same as the arithmetical
consequences of $\mathsf{ACA}$, the subsystem of second order arithmetic
obtained by adding the full scheme of induction over natural numbers (in the
language of second order arithmetic) to the well-known subsystem $\mathsf{ACA%
}_{0}$ of second order arithmetic.} Indeed, it is straightforward to
demonstrate the consistency of $\mathsf{PA}$ within the subsystem $\mathsf{CT%
}_{1}\mathsf{[PA]}$ of $\mathsf{CT[PA]}$, where $\mathsf{CT}_{n}\mathsf{[PA]}
$ is the subtheory of $\mathsf{CT[PA]}$ with the scheme of induction over
natural numbers limited to $\mathcal{L}_{\mathsf{A+T}}$-formulae that are at
most of complexity $\Sigma _{n}$ \cite[Thm.~2.8]{Lelyk+Wcislo (2017)}.
\medskip

The discussion above leaves open whether $\mathsf{CT}_{0}\mathsf{[PA]}$ is
conservative over $\mathsf{PA}$. Kotlarski \cite{Hentry 1986} established
that $\mathsf{CT}_{0}\mathsf{[PA]}$ is a subtheory of $\mathsf{CT}^{-}%
\mathsf{[PA]+Ref(PA)}$, where $\mathsf{Ref(PA)}$ is the $\mathcal{L}_{%
\mathsf{A+T}}$-sentence stating that \textquotedblleft every first order
consequence of $\mathsf{PA}$ is true\textquotedblright . Recently \L e\l yk
\cite{Lelyk Thesis} demonstrated that the converse also holds, which
immediately implies that $\mathsf{CT}_{0}\mathsf{[PA]}$ is not conservative
over \textsf{PA} since $\mathsf{Con(PA)}$ is readily provable in $\mathsf{CT}%
^{-}\mathsf{[PA]+Ref(PA).}$\footnote{%
This result was first claimed by Kotlarski \cite{Hentry 1986}, but his
proof outline of $\mathsf{Ref(PA)}$ within $\mathsf{CT}_{0}\mathsf{[PA]}$
was found to contain a serious gap in 2011 by Heck and Visser; this gap cast
doubt over the veracity of Kotlarski's claim until the issue was
resolved by \L e\l yk  in his doctoral dissertation \cite{Lelyk Thesis}. Mention should also be made of
the \L e\l yk-Wcis\l o paper \cite{Lelyk+Wcislo (2017)} which gives an
elegant proof of the non conservativity of $\mathsf{CT}_{0}\mathsf{[PA]}$
over $\mathsf{PA}$.} Kotlarski's aforementioned theorem was refined by Cie%
\'{s}li\'{n}ski \cite{Cezary 2010a} who proved that $\mathsf{CT}^{-}\mathsf{%
[PA]+}$ \textquotedblleft $\mathsf{T}$ is closed under propositional
proofs\textquotedblright\ and $\mathsf{CT}^{-}\mathsf{[PA]+Ref(PA)}$
axiomatize the same theory. The main result of this paper, in turn, refines
Cie\'{s}li\'{n}ski's result by demonstrating:\medskip

\noindent \textbf{1.~Theorem.~}$\mathsf{CT}^{-}\mathsf{[I\Delta }_{0}+%
\mathsf{Exp]+DC}$\textit{\ and }$\mathsf{CT}_{0}\mathsf{[PA]}$\textit{\
axiomatize the same first order theory.}\medskip

In the above theorem, $\mathsf{CT}^{-}\mathsf{[\mathsf{I\Delta }_{0}+\mathsf{%
Exp}]}$ is the weakening of $\mathsf{CT}^{-}\mathsf{[PA]}$ obtained by
replacing the `base theory' $\mathsf{PA}$ with its fragment consisting of
Robinson's arithmetic $\mathsf{Q}$, along with the scheme for $\Delta _{0}$%
-induction and the totality of the exponential function; and $\mathsf{DC}$
(disjunctive correctness) is the statement asserting that a disjunction of
arithmetical sentences of arbitrary finite length is true (in the sense of $%
\mathsf{T}$) iff one of the disjuncts is true. Coupled with \L e\l yk's
aforementioned result \cite{Lelyk Thesis}, Theorem 1 shows that $\mathsf{CT}%
^{-}\mathsf{[\mathsf{I\Delta }_{0}+\mathsf{Exp}]+DC,}$ $\mathsf{CT}^{-}%
\mathsf{[PA]+DC}$, and $\mathsf{CT}_{0}\mathsf{[PA]+Ref(PA)}$ are
axiomatizations of the same theory. \medskip

The plan of the paper is as follows: in Section 2 we review preliminary
definitions and results, including more precise versions of those
definitions and results mentioned in this introduction. In Section 3 we
establish the veracity of the principle $\mathsf{IC}$ (Inductive
Correctness, often referred to in the literature as \textquotedblleft
internal induction\textquotedblright ) within $\mathsf{CT}^{-}\mathsf{[%
\mathsf{I\Delta }_{0}+\mathsf{Exp}]+DC.}$ This is demonstrated by first
establishing a new general form of Visser's theorem \cite{Visser Paradox} on
nonexistence of infinite descending chains of truth definitions with the
help of (L\"{o}b's version of) G\"{o}del's second incompleteness theorem
instead of the Visser-Yablo paradox\textit{. }In Section 4 we show that $%
\mathsf{CT}_{0}\mathsf{[PA]}$ is a subtheory of $\mathsf{CT}^{-}\mathsf{%
[PA]+DC+IC}$; thus Sections 3 and 4 together constitute the proof of the
hard direction of Theorem 1 since it is routine to verify that both $\mathsf{%
DC}$ and $\mathsf{IC}$ are theorems of $\mathsf{CT}_{0}\mathsf{[PA]}$. We
close the paper with some open problems in Section 5. \medskip

\textbf{Historical Note.~}The concept of disjunctive correctness first
appeared in the work of Krajewski \cite[p.133]{Krajewski}, who called it
\textquotedblleft $\vee $-completeness\textquotedblright ; the current
terminology was coined in a working paper of Enayat and Visser that was
privately circulated in 2011, only a fragment \cite{Ali-Albert} of which has
been published so far. The working paper included the claim that $\mathsf{CT}%
^{-}\mathsf{[PA]+DC}$ is conservative over $\mathsf{PA}$, but the proof
outline presented in the paper was found in 2013 to contain a significant
gap by Cie\'{s}li\'{n}ski and his (then) doctoral students\textbf{~}\L e\l %
yk and Wcis\l o. On the other hand, in 2012 Enayat found a proof of $\mathsf{%
CT}_{0}\mathsf{[PA]}$ within $\mathsf{CT}^{-}\mathsf{[\mathsf{I\Delta }_{0}+%
\mathsf{Exp}]+DC+IC}$; his proof was only privately circulated, and later
was presented in the doctoral dissertation of \L e\l yk \cite{Lelyk Thesis}.
This proof forms the content of Section 4 of this paper. In light of these
developments, and the well-known conservativity of $\mathsf{CT}^{-}\mathsf{%
[PA]+IC}$ over $\mathsf{PA}$ (see Theorem 2.3), the question of
conservativity of $\mathsf{CT}^{-}\mathsf{[PA]+DC}$ over $\mathsf{PA}$ came
to prominence amongst truth theory experts \cite[p.226]{Cezary book}, and
has been unsuccessfully attacked by a number of researchers since 2013,
until Pakhomov established $\mathsf{IC}$ within $\mathsf{CT}^{-}\mathsf{[%
\mathsf{I\Delta }_{0}+\mathsf{Exp}]+DC}$ as in Section 3 of this paper%
\textsf{,} which, coupled with Enayat's aforementioned result, yields
Theorem 1 and exhibits the unexpected arithmetical strength of $\mathsf{DC}$%
. \bigskip

\begin{center}
\textbf{2.~PRELIMINARIES }\bigskip
\end{center}

\noindent \textbf{2.1.~Definition.}

\noindent \textbf{(a)} $\mathcal{L}_{\mathsf{A}}$ is the usual language of
first order arithmetic $\{+,\cdot ,\mathsf{S}(x),<,0\}.$ To simplify
matters, we will assume that the logical constants of first order logic
consist only of $\lnot $ (negation)$,\vee $ (disjunction)$,$ and $\exists $
(existential quantification); in particular $\forall $ (universal
quantification) as well as $\wedge $ (conjunction) and other propositional
connectives are treated here as derived notions. \medskip

\noindent \textbf{(b)} Given a language $\mathcal{L}\supseteq \mathcal{L}_{%
\mathsf{A}}$, an $\mathcal{L}$-formula $\varphi $ is said to be a $\Delta
_{0}(\mathcal{L})$\textit{-formula} if all the quantifiers of $\varphi $ are
bounded by $\mathcal{L}$-terms, i.e., they are of the form $\exists x\leq t,$
or of the form $\forall x\leq t,$ where $t$ is a term of $\mathcal{L}$ not
involving $x$. Given a predicate $\mathsf{U}(x)$, $\mathcal{L}_{\mathsf{A+U}%
} $ is the language $\mathcal{L}_{\mathsf{A}}\cup \{\mathsf{U}(x)\}$.
\medskip

\noindent \textbf{(c)} Given a language $\mathcal{L}\supseteq \mathcal{L}_{%
\mathsf{A}}$, $\mathsf{I}\Delta _{0}(\mathcal{L})$ is the scheme of
induction over natural numbers for $\Delta _{0}(\mathcal{L})$-formulae.
\textit{We shall omit the reference to }$\mathcal{L}$\textit{\ if }$\mathcal{%
L}=\mathcal{L}_{\mathsf{A}}$\textit{, e.g., a }$\Delta _{0}$\textit{-formula
is a }$\Delta _{0}(\mathcal{L}_{\mathsf{A}})$\textit{-formula; and we shall
use }$\mathsf{I}\Delta _{0}(\mathsf{U})$ \textit{to abbreviate} $\mathsf{I}%
\Delta _{0}(\mathcal{L}_{\mathsf{A+U}}).$\medskip

\noindent \textbf{(d)} $\mathsf{I}\Delta _{0}+\mathsf{Exp}$ is the fragment
of Peano arithmetic obtained by strengthening Robinson's arithmetic $\mathsf{%
Q}$ with $\mathsf{I}\Delta _{0}$ and with the sentence $\mathsf{Exp}$
that expresses the totality of the exponential function $y = 2^{x}.$ It is
well-known that $\mathsf{Exp}$ can be written as $\forall x\exists y\mathsf{%
Exp}(x,y)$, where $\mathsf{Exp}(x,y)$ is a $\Delta _{0}$-predicate which,
provably in $\mathsf{I}\Delta _{0}$, satisfies the familiar algebraic laws
governing the graph of the exponential function, cf.~\cite[Sec.~V3(c)]{Hajek-Pudlak}.\medskip

\noindent \textbf{(e)} $\mathsf{Sent}_{\mathsf{A}}(x)$ is the $\mathcal{L}_{%
\mathsf{A}}$-formula that expresses \textquotedblleft $x$ is the G\"{o}%
del-number of a formula of $\mathcal{L}_{\mathsf{A}}$ with no free
variables\textquotedblright , and $\mathsf{Form}_{\mathsf{A}}^{n}(x)$ is the
$\mathcal{L}_{\mathsf{A}}$-formula that expresses \textquotedblleft $x$ is
the G\"{o}del-number of a formula of $\mathcal{L}_{\mathsf{A}}$ with
precisely $n$ free variables\textquotedblright . We use $\mathsf{Sent}_{%
\mathsf{A}}$ and $\mathsf{Form}_{\mathsf{A}}^{n}$ to refer to the
corresponding definable classes of G\"{o}del-numbers of $\mathcal{L}_{%
\mathsf{A}}$-formulae.\medskip

\noindent \textbf{(f)} Given a (base) theory $\mathsf{B}$ whose language is $%
\mathcal{L}_{\mathsf{A}}$ and which extends $\mathsf{I}\Delta _{0}+\mathsf{%
Exp}$, $\mathsf{CT}^{-}[\mathsf{B}]$ is the theory obtained by strengthening
$\mathsf{B}$ by adding the sentences $\mathsf{tarski}_{0}$ through $\mathsf{%
tarski}_{4}$ described below, where we use the following conventions: $\tau
_{1}$ and $\tau _{2}$ vary over G\"{o}del-numbers of closed $\mathcal{L}_{%
\mathsf{A}}$-terms, $\tau _{i}^{\circ }$ is the value of the term coded by $%
\tau _{i}$, $\varphi $ and $\psi $ range over G\"{o}del-numbers of $\mathcal{%
L}_{\mathsf{A}}$-sentences, $v$ ranges over variables, $\gamma (v)$ ranges
over $\mathsf{Form}_{\mathsf{A}}^{1}$, and $\underline{x}$ is the numeral
representing the value of $x$.\medskip

\noindent $\mathsf{tarski}_{0}:=\forall x\left( \mathsf{T}(x)\rightarrow
\mathsf{Sent}_{\mathsf{A}}(x)\right) .$

\noindent $\mathsf{tarski}_{1}:=\forall \tau _{1},\tau _{1}~\left( \mathsf{T}%
(\tau _{1}=\tau _{2})\leftrightarrow \tau _{1}^{\circ }=\tau _{2}^{\circ
}\right) .$

\noindent $\mathsf{tarski}_{2}:=\forall \varphi (\left( \mathsf{T}(\lnot
\varphi )\leftrightarrow \lnot \mathsf{T}(\varphi )\right) .$

\noindent $\mathsf{tarski}_{3}:=\forall \varphi ,\psi ~\mathsf{T}(\varphi
\vee \psi )\leftrightarrow \left( \mathsf{T}(\varphi )\vee \mathsf{T}(\psi
)\right) .$

\noindent $\mathsf{tarski}_{4}:=\forall v~\forall \gamma (v)~(\mathsf{T}%
\left( \exists v~\gamma (v)\right) \leftrightarrow \exists x~\mathsf{T}%
(\gamma (\underline{x})).$\textbf{\medskip }

\noindent \textbf{(g)} $\mathsf{CT}_{0}[\mathsf{B}]:=\mathsf{CT}^{-}[\mathsf{%
B}]\cup \mathsf{I}\Delta _{0}(\mathsf{T}).$\medskip

\noindent \textbf{(h)} $\mathsf{DC}$ (disjunctive correctness) is the $%
\mathcal{L}_{\mathsf{A+T}}$-sentence asserting that $\mathsf{T}$ commutes
with disjunctions of arbitrary length, i.e., $\mathsf{DC}$ asserts that for
all numbers $s$ and for all sequences $\left\langle \varphi
_{i}:i<s\right\rangle $ from $\mathsf{Sent}_{\mathsf{A}},$ the following
equivalence holds:

\begin{center}
$\mathsf{T}(\bigvee\limits_{i<s}\varphi _{i})\leftrightarrow \exists i<s\
\mathsf{T}(\varphi _{i}),$
\end{center}

\noindent where for definiteness $\bigvee\limits_{i<s}\varphi _{i}$ is
defined$\footnote{%
One can also formulate disjunctive correctness in a stronger way by
asserting that all disjunctions (and not just the ones that are grouped to
the left) are well-behaved with respect to $\mathsf{T}$. But the current
frugal form, as shown by Theorem 1, ends up implying the seemingly stronger
form since it is easy to show in $\mathsf{CT}_{0}$ that the two forms are
equivalent.}$ by recursion: $\bigvee\limits_{i<0}\varphi _{i}:=\varphi _{0}$
and $\bigvee\limits_{i<t+1}\varphi _{i}:=\left( \bigvee\limits_{i<t}\varphi
_{i}\right) \vee \varphi _{t}$.

\noindent \textbf{(i)} We will employ the abbreviation $\bigwedge%
\limits_{i<s}\varphi _{i}$ for $\lnot \bigvee\limits_{i<s}\lnot \varphi _{i}$%
, and $\mathsf{CC}$ (conjunctive correctness) for the $\mathcal{L}_{\mathsf{%
A+T}}$-sentence

\begin{center}
$\mathsf{T}(\bigwedge\limits_{i<s}\varphi _{i})\leftrightarrow \forall i<s\
\mathsf{T}(\varphi _{i})$.
\end{center}

\begin{itemize}
\item Note that the commutativity of $\mathsf{T}$ with negation implies that
$\mathsf{DC}$ and $\mathsf{CC}$ are equivalent. \medskip
\end{itemize}

\noindent \textbf{(j)} $\mathsf{IC}$ (inductive correctness\footnote{%
This condition has been referred to as $\mathsf{Int}$ (internal induction)
in the literature.}) is the $\mathcal{L}_{\mathsf{A+T}}$-sentence asserting
that each $\mathcal{L}_{\mathsf{A}}$-instance of induction over natural
numbers is true, i.e., $\mathsf{IC}$ asserts that for all unary $\mathcal{L}%
_{\mathsf{A}\text{ }}$-formulae $\psi =\psi (x),$ $\mathsf{T}(\ulcorner
\mathsf{Ind}_{\psi }\urcorner )$ holds, where $\mathsf{Ind}_{\psi }$ is the
following $\mathcal{L}_{\mathsf{A}}$-sentence that asserts that $\psi $ is
inductive:

\begin{center}
$\psi (\underline{0})\rightarrow \left( \forall x\left( \psi (x)\rightarrow
\psi (x+1)\right) \rightarrow \forall x\ \psi (x)\right) .$
\end{center}

The $\mathsf{B}=\mathsf{PA}$ case of Theorem 2.2 below, and its elaboration
Theorem 2.3, were first established in the work of Krajewski, Kotlarski, and
Lachlan \cite{KKL} for $\mathsf{B}=\mathsf{PA}$, where $\mathsf{PA}$ is
formulated in a relational language, and `domain constants' play the role of
numerals. As mentioned in the introduction to this paper, their result was
couched in model theoretic terms involving the notions of recursive
saturation and satisfaction classes, but it is well-known that their
formulation is equivalent to appropriately formulated conservatity
assertions (the key ingredients of this equivalence are the following facts:
Every consistent theory in a countable language has a recursively saturated
model, and countable recursively saturated models are resplendent). Later
Kaye \cite{Kaye book} developed the theory of satisfaction classes over
models of $\mathsf{PA}$ in languages incorporating function symbols; his
work was extended by Engstr\"{o}m \cite{Engstrom} to truth classes over
models of $\mathsf{PA}$ in functional languages.\footnote{%
The subtle distinction between satisfaction classes and truth classes, and
their close relationship, is explained in \cite{Ali-Albert} and \cite[Ch.~7]%
{Cezary book}.} More recently, newer and more informative proofs of Theorems
2.2 and 2.3 have been found in the joint work of Visser and Enayat \cite%
{Ali-Albert} (with base theories that support a modicum of coding, and which
are formulated in a relational language), and by Leigh \cite{Leigh} (for
functional base theories that support a modicum of coding). As verified by
Cie\'{s}li\'{n}ski \cite[Ch.~7]{Cezary book}, the Visser-Enayat model
theoretic methodology can be extended so as to accommodate functional
languages. \medskip

\noindent \textbf{2.2.~Theorem.~}$\mathsf{CT}^{-}[\mathsf{B]}$ \textit{is
conservative over }$\mathsf{B}$ \textit{for every arithmetical base theory }$%
\mathsf{B}$\textit{\ extending} $\mathsf{I}\Delta _{0}+\mathsf{Exp}$.\medskip

\noindent \textbf{2.3.~Theorem.~}$\mathsf{CT}^{-}[\mathsf{PA]+IC}$ \textit{%
is conservative over} $\mathsf{PA.}$\medskip

The direction $(a)\Rightarrow (b)$ of Theorem 2.4 below is due to Kotlarski
\cite{Hentry 1986}; the other direction is due to \L e\l yk \cite{Lelyk
Thesis}.\medskip

\noindent \textbf{2.4.~Theorem}.\textbf{~}(Kotlarski-\L e\l yk) \textit{The
following theories are deductively equivalent}:\medskip

\noindent \textbf{(a) }$\mathsf{CT}^{-}[\mathsf{PA]}+$\textsf{\ }$\mathsf{%
Ref(PA).\medskip }$

\noindent \textbf{(b)} $\mathsf{CT}_{0}[\mathsf{PA].}$\medskip

The direction $(a)\Rightarrow (b)$ of Theorem 2.5 below is due to Cie\'{s}li%
\'{n}ski \cite{Cezary 2010a}, who refined Kotlarski's proof of the direction
$(a)\Rightarrow (b)$ of Theorem 2.4; the other direction involves a routine
induction.\medskip

\noindent \textbf{2.5.}~\textbf{Theorem.~}(Cie\'{s}li\'{n}ski) \textit{The
following theories are deductively equivalent}:\medskip

\noindent \textbf{(a)} $\mathsf{CT}^{-}[\mathsf{PA]}+``\mathsf{T}$ is closed
under propositional proofs\textquotedblright .$\mathsf{\medskip }$

\noindent \textbf{(b)} $\mathsf{CT}_{0}[\mathsf{PA].}$\bigskip

\begin{center}
\textbf{3.~DISJUNCTIVE CORRECTNESS IMPLIES INDUCTIVE CORRECTNESS\ }\bigskip
\end{center}

\noindent \textbf{3.1.}~\textbf{Definition.~}$\mathsf{ITB}$ (iterated truth
biconditionals) is a theory formulated in two-sorted first order logic. The
first sort $x,y,z\ldots $ of $\mathsf{ITB}$ is for the `natural numbers'.
The second sort $\alpha ,\beta ,\gamma ,\ldots $ is for the \textit{indices}
of truth definition. The language $\mathcal{L}_{\mathsf{ITB}}$ of $\mathsf{%
ITB}$ is obtained by augmenting the language $\mathcal{L}_{\mathsf{A}}$ of
arithmetic with two binary predicates: $\alpha \prec \beta $ and $\mathsf{T}%
(\alpha ,x)$, but we shall write $\mathsf{T}(\alpha ,x)$ as $\mathsf{T}%
_{\alpha }(x)$ to display the indexicality of $\alpha $. The axioms of
\textsf{ITB}\ come in three groups. The first group consists of the axioms
of $\mathsf{Q}$ (Robinson arithmetic); the second group stipulates that $%
\prec $ is a strict partial order; and the third group consists of the
following biconditionals $\mathsf{B}_{\varphi }$:
\begin{equation*}
\mathsf{B}_{\varphi }:=\forall \alpha (\mathsf{T}_{\alpha }(\gnmb{\varphi})%
\mathrel{\leftrightarrow}\mathbb{\varphi }^{\prec \alpha }),
\end{equation*}%
where $\mathbb{\varphi }$ ranges over all $\mathcal{L}_{\mathsf{ITB}}$%
-sentences, and for each index variable $\alpha $, $\mathbb{\varphi }^{\prec
\alpha }$ denotes the relativization of $\mathbb{\varphi }$ to the cone of
indices below $\alpha $, i.e. the formula obtained by replacing all the
quantifiers of the form $\forall \beta $ $(\exists \beta )$ with $\forall
\beta \prec \alpha $ $(\exists \beta \prec \alpha )$, and if there is a
bounded instance of $\alpha $ we make the appropriate renaming. Clearly $%
\mathbb{\varphi }^{\prec \alpha }=\varphi $ if $\varphi $ is a purely
arithmetical formula. \medskip

\begin{itemize}
\item Note that we take the theory $\mathsf{ITB}$ over the variant of
many-sorted logic that allows domains of some sorts to be empty.

\item We will use the following convention to lighten the notation: The notation $\ulcorner \varphi \urcorner $ for the G\"{o}del
number of a formula $\varphi $ will be generally used, but the corner-notation will be omitted when $%
\varphi $ appears inside of a truth predicate $\mathsf{T}$, or inside an
indexed version of $\mathsf{T}.$
\end{itemize}

The proof of the following theorem was inspired by the recent James Walsh
proof \cite{Pakhomov-Walsh} of nonexistence of infinite recursive provably
descending chains of sentences with respect to $<_{\mathsf{Con}}$%
-order.\medskip

\noindent \textbf{3.2.~Theorem.~}$\mathsf{ITB}+\exists \alpha (\alpha
=\alpha )$ \textit{proves the existence of a }$\prec $\textit{-minimal
element}. \textit{Equivalently, the following theory }$\mathsf{DTB}$ (%
\textit{descending truth biconditionals}) \textit{is inconsistent:}

\begin{center}
$\mathsf{DTB:=ITB}+\forall \alpha \exists \beta (\beta \prec \alpha
)+\exists \alpha (\alpha =\alpha )$\textit{.} \medskip
\end{center}

\noindent \textbf{Proof.~}Consider the formula $\mathsf{\theta }(x)$:
\begin{equation*}
\mathsf{\theta }(x):=\forall \alpha (\mathsf{T}_{\alpha }(x)).
\end{equation*}%
\noindent We will verify that $\mathsf{\theta }(x)$ satisfies the following
Hilbert-Bernays-L\"{o}b derivability conditions in $\mathsf{DTB}$; in what
follows $\mathsf{\varphi }$ and $\mathsf{\psi }$ range over all sentences of
the language of $\mathsf{DTB}$:\medskip

\noindent HBL-1. $\mathsf{DTB}\vdash \varphi \Longrightarrow \mathsf{DTB}%
\vdash \mathsf{\theta }(\ulcorner \varphi \urcorner )$.

\noindent HBL-2. $\mathsf{DTB}\vdash \mathsf{\theta }(\ulcorner \varphi
\rightarrow \psi \urcorner )\rightarrow (\mathsf{\theta }(\ulcorner \varphi
\urcorner )\rightarrow \mathsf{\theta }(\ulcorner \psi \urcorner ))$.

\noindent HBL-3. $\mathsf{DTB}\vdash \mathsf{\theta }(\ulcorner \varphi
\urcorner )\rightarrow \mathsf{\theta }(\ulcorner \mathsf{\theta (}\ulcorner
\varphi \urcorner )\urcorner )$.\medskip

\noindent It is easy to see that for every axiom $\mathsf{A}$ of $\mathsf{DTB%
}$ we have $\mathsf{DTB}\vdash \mathsf{A}^{\prec \alpha }$. The assumption
embodied in $\mathsf{DTB}$ that the set of indices is nonempty and has no
minimal element is invoked only at this step of the proof, in verifying that
$\mathsf{DTB}\vdash \mathsf{A}^{\prec \alpha }$ for

\begin{center}
$\mathsf{A:}=\forall \gamma \exists \beta (\beta \prec \gamma )\wedge
\exists \gamma (\gamma =\gamma ).$
\end{center}

\noindent Routine considerations show that HBL-1 holds based on an easy
induction on the number of steps of the proof of $\varphi $ within $\mathsf{%
DTB.}$\medskip

\noindent For a given $\varphi $ and $\psi ,$ HBL-2 follows directly from
the biconditional axioms $\mathsf{B}_{\varphi \rightarrow \psi }$, $\mathsf{B%
}_{\varphi },$ and $\mathsf{B}_{\psi }$ of $\mathsf{ITB.}$\medskip

\noindent Finally, HBL-3 holds since:

\begin{center}
$\mathsf{ITB}\vdash \mathsf{\theta }(\ulcorner \mathsf{\theta (}\ulcorner
\varphi \urcorner )\urcorner )\longleftrightarrow \left( \forall \alpha
\forall \beta ((\beta \prec \alpha )\rightarrow \mathsf{T}_{\beta }(\varphi
))\right) $.
\end{center}

\noindent On the other hand, the formula $\lnot \theta (\ulcorner \underline{%
0}=\underline{1}\urcorner )$ is provable in $\mathsf{ITB}$, and therefore in
$\mathsf{DTB}$. So by L\"{o}b's version\footnote{%
L\"{o}b's paper \cite{Lob}, in which the venerable `L\"{o}b's Theorem' was
proved, is responsible for the now common standard textbook framework for
the presentation of `abstract' form of G\"{o}del's second incompleteness
theorem: If $T$ is a theory extending $\mathsf{Q}$ that supports a unary
predicate $\theta (x)$ satisfiying condition HBL-1, HBL2, and HBL-3, and $%
\gamma $ is a $T$-provable fixed point of $\lnot \theta $ (i.e., the
equivalence of $\gamma $ and $\lnot \theta (\ulcorner \gamma \urcorner )$ is
provable in $T$), then either $\lnot \exists x\ \theta (x)$ is not provable
in $T$, or $T$ is inconsistent. See, e.g., \cite[Ch.~18]{Boolos et al.}, for
the presentation of such a general form of G\"{o}del's second incompleteness
theorem.} of G\"{o}del's second incompleteness theorem, $\mathsf{DTB}$ is
inconsistent. \hfill $\square $\medskip

\noindent \textbf{3.3.~Lemma}.\textbf{~}$\mathsf{CT}^{-}\mathsf{[\mathsf{I}%
\Delta _{0}+\mathsf{Exp}]}+\mathsf{DC}$ \textit{proves} $\mathsf{IC}$%
.\medskip

\noindent \textbf{Proof.~}By Theorem 3.2 we can fix an inconsistent finite
subtheory $\mathsf{DTB}^{-}$ of $\mathsf{DTB}$. Suppose $\mathsf{DTB}^{-}$
contains only biconditionals for the formulae $\mathsf{\varphi }_{0},\ldots ,%
\mathsf{\varphi }_{k-1}$. We will use $\mathsf{ITB}^{-}$ to denote the
subtheory of $\mathsf{ITB}$ whose only biconditional axioms are $\left\{
\mathsf{B}_{\varphi _{i}}:i<k\right\} .$\medskip

For the rest of the proof we will reason in $\mathsf{CT}^{-}[\mathsf{I}%
\Delta _{0}+\mathsf{Exp}]+\mathsf{DC}$. In order to prove $\mathsf{IC}$ we
assume for a contradiction that some arithmetical $\mathsf{\psi }(x)$ is not
inductive in the sense of $\mathsf{T}$, i.e., we have:

\begin{center}
$\lnot \mathsf{T}\left( \mathsf{\psi }(0)\rightarrow \left( \forall x\left(
\mathsf{\psi }(x)\rightarrow \mathsf{\psi }(x+1)\right) \rightarrow \forall
x\;\mathsf{\psi }(x)\right) \right) $.
\end{center}

\noindent By induction on natural numbers $n$ we define interpretations $%
\iota _{n}$ of $\mathsf{ITB}^{-}$. Note that each interpretation is just a
finite sequence of formulae and thus could be easily represented in
arithmetic. The interpretation of arithmetic in each $\iota _{n}$ is the
identity interpretation, but the domain of indices of truth definition $%
\iota _{n}$ is given by the formula $D^{(n)}(x)$:
\begin{equation*}
x<\underline{n}\wedge \lnot \psi (x).
\end{equation*}%
For all $n$ the relation $\prec $ is interpreted by $<$. The formula $%
\mathsf{T}_{\alpha }(x)$ is interpreted by the formula $\mathsf{IT}%
^{(n)}(y,x)$, where $y$ corresponds to $\alpha $, and $x$ corresponds to
itself:
\begin{equation*}
\bigwedge\limits_{i<k}(x=\gnmb{\varphi_i}\rightarrow
\bigwedge\limits_{m<n}((y=\underline{m}\wedge \lnot \psi (\underline{m}%
))\rightarrow \iota _{m}(\mathsf{\varphi }_{i}))),
\end{equation*}%
where $\iota _{m}(\mathsf{\varphi }_{i})$ is the $\iota _{m}$-translation of
the sentence $\mathsf{\varphi }_{i}$. It is easy to see that this definition
could be carried out in $\mathsf{\mathsf{I\Delta }_{0}+\mathsf{Exp}}$%
.\medskip

Let us now prove that the translations given by $\iota _{n}$ are indeed the
desired interpretations inside $\mathsf{T}$, i.e., we need to prove that for
all $n$ and axioms $\mathsf{A}$ of $\mathsf{ITB}^{-}$ we have $\mathsf{T}%
(\iota _{n}(\mathsf{A}))$. Clearly it is the case for all the axioms of $%
\mathsf{Q}$ and the axioms of partial order for $\prec $. Now let us show
that for any $s<k:$

\begin{description}
\item[$\mathbf{(\blacktriangle )}$] $\quad \mathsf{T}(\iota_n(\forall \alpha
(\mathsf{T}_{\alpha}(\mathsf{\varphi}_s)\mathrel{\leftrightarrow} \mathsf{%
\varphi}_s^{\prec \alpha}))).$
\end{description}

By compositional axioms, we just need to show that for all $u$ such that $%
u<n $ and $\mathsf{T}(\lnot \psi (\underline{u}))$ we have:
\begin{equation*}
\mathsf{T}\left( \bigwedge\limits_{i<k}\left( \mathsf{\varphi }_{s}=\mathsf{%
\varphi }_{i}\rightarrow \bigwedge\limits_{m<n}\left( \left( \underline{u}=%
\underline{m}\wedge \lnot \mathsf{\psi }(\underline{m})\right) \rightarrow
\iota _{m}(\mathsf{\varphi }_{i})\right) \right) \right) \mathrel{%
\leftrightarrow}\mathsf{T}\left( \iota _{n}\left( \mathsf{\varphi }%
_{s}^{\prec \underline{u}}\right) \right) .
\end{equation*}%
Now by compositional axioms and $\mathsf{DC}$ (in the form of $\mathsf{CC}$, as explained in part (i) of Definition 2.1)
our task can be reduced to proving the equivalence:%
\begin{equation*}
\mathsf{T}\left( \iota _{u}\left( \mathsf{\varphi }_{s}\right) \right) %
\mathrel{\leftrightarrow}\mathsf{T}\left( \iota _{n}\left( \mathsf{\varphi }%
_{s}^{\prec \underline{u}}\right) \right) .
\end{equation*}%
In order to prove this we will show by induction on subformulae $\mathsf{%
\theta }$ of $\mathsf{\varphi }_{s}$ that for the universal closure $%
\overline{\mathsf{\theta }}$ of $\mathsf{\theta }$:
\begin{equation*}
\mathsf{T}\left( \iota _{u}(\overline{\mathsf{\theta }}\right) %
\mathrel{\leftrightarrow}\mathsf{T}\left( \iota _{n}(\overline{\mathsf{%
\theta }}^{\prec \underline{u}})\right) .
\end{equation*}%
Note that since $\mathsf{\varphi }_{s}$ is a fixed formula with finitely
many subformulae, actually this external induction will be formalizable in $%
\mathsf{CT}^{-}[\mathsf{I\Delta }_{0}+\mathsf{Exp}]+\mathsf{DC}$ despite the
fact that it lacks the induction axiom for the appropriate class of
formulae. The only non-trivial case here is the case when $\mathsf{\theta }$
is $\mathsf{T}_{\alpha }(x)$:%
\begin{equation*}
\mathsf{T}\left( \iota _{u}\left( \forall \alpha \forall x\mathsf{T}_{\alpha
}(x)\right) \right) \mathrel{\leftrightarrow}\mathsf{T}(\iota _{n}(\forall
\alpha \prec \underline{u}~\forall x(\mathsf{T}_{\alpha }(x)))).
\end{equation*}%
Hence we just need to show that for all $p<u$ such that $\mathsf{T}(\lnot \mathsf{%
\psi }(\underline{p}))$, and for all $o$, the following pair of formulae (whose only
formal difference is in the bound for indices of the second conjunction) are
equivalent:
\begin{equation*}
\mathsf{T}\left( \bigwedge\limits_{i<k}\left( \underline{o}=\mathsf{\varphi }%
_{i}\rightarrow \bigwedge\limits_{m<u}\left( \left( \underline{p}=\underline{%
m}\wedge \lnot \mathsf{\psi }(\underline{m})\right) \rightarrow \iota _{m}(%
\mathsf{\varphi }_{i})\right) \right) \right) ,
\end{equation*}%
\begin{equation*}
\mathsf{T}\left( \bigwedge\limits_{i<k}\left( \underline{o}=\mathsf{\varphi }%
_{i}\rightarrow \bigwedge\limits_{m<n}\left( \left( \underline{p}=\underline{%
m}\wedge \lnot \mathsf{\psi }(\underline{m})\right) \rightarrow \iota _{m}(%
\mathsf{\varphi }_{i})\right) \right) \right) .
\end{equation*}%
But since $p<u<n$, we trivially use $\mathsf{DC}$ (in the form of $\mathsf{CC%
}$) to show that the formulae
are indeed equivalent. Thus we conclude that $(\blacktriangle )$
holds.\medskip

Choose some $n$ such that $\mathsf{T}(\lnot \psi(\underline{n}))$; this is
possible since we assumed that induction fails for $\psi (x)$ in the sense
of $\mathsf{T}$. It is easy to see that $\iota _{n}$ actually is an
interpretation of $\mathsf{DTB}^{-}$ inside $\mathsf{T}$. Now we could
follow the proof of inconsistency of $\mathsf{DTB}^{-}$ inside $\mathsf{T}$
to derive a contradiction, thereby completing the proof of $\mathsf{IC}$%
.\hfill $\square \medskip $

\noindent \textbf{3.4.~Corollary}.\textbf{~}$\mathsf{CT}^{-}\mathsf{[\mathsf{%
I}\Delta _{0}+\mathsf{Exp}]}+\mathsf{DC}$ \textit{proves} $\mathsf{PA}$,
\textit{and therefore} $\mathsf{CT}^{-}\mathsf{[\mathsf{I}\Delta _{0}+%
\mathsf{Exp}]}+\mathsf{DC}$ \textit{and} $\mathsf{CT}^{-}\mathsf{[PA]}+%
\mathsf{DC}$ \textit{axiomatize the same theory}.\medskip

\noindent \textbf{Proof.~}This is an immediate consequence of Lemma 3.3, and
the provability of Tarski bi-conditionals in $\mathsf{CT}^{-}\mathsf{[%
\mathsf{I}\Delta _{0}+\mathsf{Exp}]}$.\hfill $\square $\bigskip

\begin{center}
\textbf{4.~DISJUNCTIVE CORRECTNESS + INDUCTIVE CORRECTNESS IMPLIES }$\Delta
_{0}(\mathsf{T})$\textbf{-INDUCTION\ }\bigskip
\end{center}

In this section we shall prove that $\mathsf{I}\Delta _{0}(\mathsf{T})$ is
provable in $\mathsf{CT}^{-}[\mathsf{PA}]+\mathsf{DC+IC},$ which, coupled
with Lemma 3.3 completes the proof of the nontrivial direction of Theorem 1.
We begin with a key definition:\medskip

\noindent \textbf{4.1.~Definition.~}In what follows $\in _{\mathsf{Ack}}$ is
\textquotedblleft Ackermann's epsilon\textquotedblright , i.e., $x\in _{%
\mathsf{Ack}}y$ is the arithmetical formula that expresses \textquotedblleft
the $x$-th bit of the binary expansion of $y$ is a 1\textquotedblright .
\medskip

\noindent \textbf{(a)} For a unary predicate $\mathsf{U}(x),$ the $\mathcal{L%
}_{\mathsf{A}+\mathsf{U}}$ sentence $\mathsf{PC}(\mathsf{U})$ (read as
\textquotedblleft $\mathsf{U}$ is piece-wise coded\textquotedblright ) is
the following sentence:

\begin{center}
$\forall u\ \exists y\ \forall x\ \left[ (\mathsf{U}(x)\wedge
x<u)\leftrightarrow x\in _{\mathsf{Ack}}y\right] .$
\end{center}

\noindent \textbf{(b)} Given an $n$-ary formula $\mathcal{L}_{\mathsf{A+U}}$%
-formula $\varphi (\mathsf{U},x_{0},...,x_{n-1}),$ $\mathsf{PC}_{\varphi }$
is the following $\mathcal{L}_{\mathsf{A+U}}$-sentence:

\begin{center}
$\forall u\ \exists y\ \forall x_{0},...\forall x_{n-1}\ \left[ (\varphi (%
\mathsf{U},x_{0},...,x_{n-1})\wedge \vec{x}<u)\leftrightarrow \vec{x}\in _{%
\mathsf{Ack}}y\right] ,$
\end{center}

\noindent where $\vec{x}=\left\langle x_{i}:i<n\right\rangle $ is a
canonical code for the ordered $n$-tuple $\left( x_{0},...,x_{n-1}\right)$.
\medskip

The following lemma shows that over $\mathsf{\mathsf{I}\Delta _{0}+\mathsf{%
Exp}}$ the scheme $\mathsf{I}\Delta _{0}(\mathsf{U})$ is equivalent to the
single sentence `$\mathsf{U}$ is piecewise coded'\textsf{.} The lemma is
folklore; we present the proof for the sake of completeness.\medskip

\noindent \textbf{4.2. Lemma.~}\textit{The following are equivalent over} $%
\mathsf{\mathsf{I}\Delta _{0}+\mathsf{Exp}}$:\medskip

\noindent $\mathbf{(a)}$\textbf{\ }$\mathsf{I}\Delta _{0}(\mathsf{U}).$%
\medskip

\noindent $\mathbf{(b)}$ $\mathsf{PC(U)}$.\medskip

\noindent \textbf{Proof. }We will reason in $\mathsf{\mathsf{I}\Delta _{0}+%
\mathsf{Exp}}$. Recall that $\in _{\mathsf{Ack}}$ has a $\Delta _{0}$%
-definition within $\mathsf{I}\Delta _{0}$ \cite{Hajek-Pudlak}.\medskip

\noindent $(\mathbf{a}\longrightarrow \mathbf{b})$: Assume $\mathsf{I}\Delta
_{0}(\mathsf{U}).$ Given $u$, let $w$ be the Ackermann-code for the set of
predecessors of $u.$ Clearly $w=\sum\limits_{i<u}2^{i}=2^{u}-1,$ and $w$ is
an upper bound for any $w^{\prime }$ that codes a subset of the predecessors
of $u$. Consider the $\Delta _{0}$-formula:

\begin{center}
$\delta (u):=\exists y<w+1\ \forall x<u\ \left[ (\mathsf{U}(x)\wedge
x<u)\leftrightarrow x\in _{\mathsf{Ack}}y\right] .$
\end{center}

\noindent A simple induction using $\mathsf{I}\Delta _{0}(\mathsf{U})$ shows
that $\forall u\ \delta (u)$ holds.\medskip

\noindent $(\mathbf{b}\rightarrow \mathbf{a})$: A straightforward induction
on the complexity of $\Delta _{0}(\mathsf{U})$-formulae shows that if $%
\mathsf{U}$ is piecewise coded and $\delta (\mathsf{U},x_{0},...,x_{n-1})$
is a $\Delta _{0}(\mathsf{U})$-formula, then $\mathsf{PC}_{\delta }$ holds.
We could then trivially deduce the least number principle for $\Delta _{0}(%
\mathsf{U})$-formula, which is of course equivalent to $\mathsf{I}\Delta
_{0}(\mathsf{U}).$\hfill $\square $\medskip

\begin{itemize}
\item In the next lemma and its proof, $\mathsf{Code}(c,\varphi ,u)$ denotes
the ternary $\mathcal{L}_{\mathsf{A}+\mathsf{T}}$-formula

$\forall x<u\left( \mathsf{T}\left( \varphi (\underline{x})\right)
\leftrightarrow x\in _{\mathsf{Ack}}c\right) ,$ and $\mathsf{PC}(\varphi )$
denotes the $\mathcal{L}_{\mathsf{A}+\mathsf{T}}$-formula $\forall u\
\exists c\ \mathsf{Code}(c,\varphi ,u).$
\end{itemize}

\noindent \textbf{4.3.~Lemma.~}$\mathsf{CT}^{-}[\mathsf{\mathsf{I}\Delta
_{0}+\mathsf{Exp}}]+\mathsf{IC}$ \textit{proves} $\forall \varphi (\mathsf{%
Form}_{\mathsf{A}}^{1}(\varphi )\rightarrow $ $\mathsf{PC}(\varphi )).$%
\medskip

\noindent \textbf{Proof. }We reason in $\mathsf{CT}^{-}[\mathsf{\mathsf{I}%
\Delta _{0}+\mathsf{Exp}}]+\mathsf{IC.}$ Given $\varphi (x)$ in $\mathsf{Form%
}_{\mathsf{A}}^{1}$, we need to show:\medskip

\noindent (1) $\forall u\ \exists c\ \mathsf{Code}(c,\varphi ,u).$\medskip

\noindent By the compositional properties of $\mathsf{T}$, (1) is equivalent
to:\medskip

\noindent (2) $\mathsf{T}(\forall u~\psi (u)),$ where $\psi (u):=\exists
c(\forall x<u~\varphi (x)\leftrightarrow x\in _{\mathsf{Ack}}c).$\medskip

\noindent On the other hand, $\forall u~\psi (u)$ is the conclusion of the
formula $\mathsf{Ind}_{\psi }$ (asserting the inductive property of $\psi $)
given by $\mathsf{IC}$. So (2) follows directly from $\mathsf{IC}$ and the
easily verified facts $\mathsf{T}(\psi (\underline{0}))$ and $\mathsf{T}%
(\forall u~\left( \psi (u)\rightarrow \psi (u+1)\right) ).$\hfill $\square $%
\medskip

\noindent \textbf{4.3.1.~Remark.~}Lemma 4.3 can be readily strengthened to a
more general result whose proof we leave to the reader: the $\omega $%
-interpretation that can be informally described by the motto
\textquotedblleft sets are $\mathsf{T}$-extensions of unary arithmetical
formulae\textquotedblright\ satisfies $\mathsf{ACA}_{0}$, provably in $%
\mathsf{CT}^{-}[\mathsf{PA}]+\mathsf{IC.}$ Indeed, it is a theorem of $%
\mathsf{CT}^{-}[\mathsf{PA}]$ that $\mathsf{IC}$ is equivalent to the
veracity of $\mathsf{ACA}_{0}$ within this interpretation. \medskip

\noindent \textbf{4.4.} \textbf{Lemma.~}$\mathsf{CT}^{-}[\mathsf{\mathsf{%
I\Delta }_{0}+\mathsf{Exp}}]+\mathsf{DC+IC}\vdash \mathsf{I}\Delta _{0}(%
\mathsf{T}).$\medskip

\noindent \textbf{Proof. }Reason in\textbf{\ }$\mathsf{CT}^{-}[\mathsf{%
\mathsf{I\Delta }_{0}+\mathsf{Exp}}]+\mathsf{DC+IC}$. By Lemma 4.2, it
suffices to show that $\mathsf{T}$ is piecewise coded. Let $\left\langle
\varphi _{i}:i<u\right\rangle $ be the sequence of arithmetical sentences
such that $\varphi _{i}$ is the sentence with G\"{o}del-number $i$ if there
is such a sentence, and otherwise $\varphi _{i}$ is the sentence $\underline{%
0}=\underline{1}$. We wish to show that $\left\{ i<u:\mathsf{T}(\varphi
_{i})\right\} $ is coded. Towards this goal, consider the unary formula $%
\theta (x)\in \mathsf{Form}_{\mathsf{A}}$ given by:
\begin{equation*}
\theta (x):=\bigvee_{i<u}\left( (x=\underline{i})\wedge \varphi _{i}\right) .
\end{equation*}

\noindent \textbf{Claim }$\mathbf{(\ast )}$ $\forall i<u~\left[ \mathsf{T}%
(\varphi _{i})\leftrightarrow \mathsf{T}(\theta (\underline{i}))\right] .$%
\medskip

\noindent $(\rightarrow )$ Suppose $\mathsf{T}(\varphi _{i})$ for some $i<u.$
Then $\mathsf{T}\left( (\underline{i}=\underline{i})\wedge \varphi
_{i}\right) ,$ and hence by $\mathsf{DC}$ we have$\mathsf{T}(\theta (%
\underline{i})).\medskip $

\noindent $(\leftarrow )$ Suppose $\mathsf{T}(\theta (\underline{i}))$ for
some $i<u.$ Then by $\mathsf{DC}$, there is some $j<u$ such that $\mathsf{T}%
\left( (\underline{i}=\underline{j})\wedge \varphi _{j}\right) .$ So $%
\mathsf{T}(\varphi _{i})$ holds since $\mathsf{T}$ commutes with conjunction
and $\mathsf{T}(\underline{i}=\underline{j})$ holds iff $i=j.$\medskip

\noindent By coupling Claim $(\ast )$ together with Lemma 4.3, we can
conclude that $\left\{ i<u:\mathsf{T}\left( \theta (\underline{i})\right)
\right\} $ is coded. \hfill $\square $\bigskip

\begin{center}
\textbf{5.~CLOSING\ REMARKS\ AND OPEN QUESTIONS}\bigskip
\end{center}

\noindent \textbf{5.1.}~\textbf{Question. }\textit{Is the generalization of
Theorem 1 in which} \textsf{CT}$^{-}$ \textit{is weakened to} $\mathsf{CS}%
^{-}$ (\textit{where} $\mathsf{S}$ \textit{stands for satisfaction}) \textit{%
true}?

\begin{itemize}
\item The notion $\mathsf{CS}^{-}[\mathsf{B}]$ is defined in \cite%
{Ali-Albert} for base theories $\mathsf{B}$ formulated in relational
languages, using the notation $\mathsf{B}^{\mathsf{FS}}$ ($\mathsf{FS}$ for
\textquotedblleft full satisfaction"); and in \cite[Ch.~7]{Cezary book}
for functional languages. We expect that an examination of the proofs in
Sections 3 and 4 would show that this question has a positive answer.
\end{itemize}

\noindent \textbf{5.2.}~\textbf{Question.~}\textit{Is} $\mathsf{IC}$ \textit{%
provable in} $\mathsf{CT}^{-}[\mathsf{S}_{2}^{1}]+\mathsf{DC}?$

\begin{itemize}
\item In the above, $\mathsf{S}_{2}^{1}$ is Buss's well-known arithmetical theory whose provable recursive functions are precisely the functions computable in polynomial time, as in \cite{Buss}. For the above question to make sense, part (f) of Definition 1.1 should
be adjusted so as to accommodate the fact that the language of $\mathsf{S}%
_{2}^{1}$ extends $\mathcal{L}_{\mathsf{A}}.$ In the proof of Lemma 3.3,
most likely it is possible to use some tricks with effective formulae in
order to modify the definition of $\iota _{n}$ in such a way that their
sizes will be polynomial. But in order for the construction to work we will
also need to ensure that $\mathsf{DC}$ is still enough to show that $\iota
_{n}$ are indeed interpretations inside the truth predicate.
\end{itemize}

\noindent \textbf{5.3.}~\textbf{Question.~} \textit{Is there a fixed-point
free proof of Lemma }3.3?

\begin{itemize}
\item The proof of Lemma 3.3 is based on Theorem 3.2, whose proof implicitly
relies (in the very last step) on the existence of a fixed point for the
formula $\lnot \theta (x).$
\end{itemize}

\noindent \textbf{5.4.~Question.~}\textit{Is} $\mathsf{CT}^{-}[\mathsf{PA}]+%
\mathsf{\ ``T}$ contains all arithmetical instances of propositional
tautologies\textquotedblright\ \textit{conservative over} $\mathsf{PA}$?

\begin{itemize}
\item Cie\'{s}li\'{n}ksi (\cite{Cezary 2010b}, \cite[Ch.~12]{Cezary book})
has shown that $\mathsf{CT}^{-}[\mathsf{PA}]+\mathsf{\ ``T}$ contains all
arithmetical instances of theorems of first order logic\textquotedblright\ is deductively equivalent
to $\mathsf{CT}^{-}[\mathsf{PA}]+\mathsf{Ref(PA)}$, and therefore by Theorem
2.4, it is yet another axiomatization of $\mathsf{CT}_{0}[\mathsf{PA}].$
\end{itemize}

\begin{center}
\bigskip
\end{center}

\end{document}